# Effects on two semi-circular wall heaters in a rectangular enclosure containing trapezoidal heated obstacle in presence of MHD

Sayeda Fahmida Ferdousi[1], Md. Abdul Alim[1]

[1] *Department of Mathematics, Bangladesh University of Engineering and Technology, Dhaka-1000, Bangladesh*

**Abstract**

A computational study has been done to illustrate the effects on two semi-circular wall heaters placed in a mixed convection rectangular enclosure containing trapezoidal heated obstacle in presence of magnetic field. The upper wall moves with a velocity in the positive $x$-direction. Trapezoidal obstacle is located in the middle of the cavity and two semi-circular wall heaters are placed at the bottom wall with fixed distance between them. These two horizontal walls are kept adiabatic except two semi-circular wall heaters. Two vertical walls are kept at low temperature and concentration and the trapezoidal obstacle and two semi-circular walls are kept at high temperature and concentration. Finite element method is used to solve governing equations. The present analysis is performed for parameters such as Hartmann number, Buoyancy ratio and Richardson number. The effects of aforesaid parameters are explored on the fluid flow and temperature at two semi-circular wall heaters. The results show that heat transfer rate along right semi-circular wall heater dominates the left semi-circular wall heater with the increase of both the parameters Hartmann number and Buoyancy ratio.

**1. Introduction:**
Mixed convection flow with heat and mass transfer in lid-driven cavities have been receiving a considerable attention due to the attention of shear flow caused by the movement of moving wall and the combined effects of thermal and species diffusion. Such problems are commonly encountered in many engineering, technological and natural processes. The problem is chosen as benchmark case for the solution of Navier-stokes equations by U. Ghial et al. (1982). The solution of the problem is important to understand the phenomena of underground heat pumps, lakes and reservoirs.
Heat and fluid flow in lid-driven enclosures with different geometry and thermal boundary conditions and lid position were studied by many researchers like Dong et al. (2004),

Hasanuzzaman et al. (2012), Islam et al. (2012) and Rahman et al (2012). In this context, Al-Salem at al. (2012) studied the effects of moving lid direction on MHD mixed convection in a linearly heated cavity. Mixed convection flows within a square cavity with linearly heated side wall was studied by Basak at al. (2009). Irregular surface has vast application in engineering problem studied by Al-Amiri et al. (2007). His study exhibits the effects of sinusoidal wavy bottom surface on mixed convection heat transfer in al id-driven cavity. Combined effect of Hartmann and Rayleigh numbers on free convective flow in a square cavity with different positions of heated elliptic obstacle was studied by Bhuiyan at al. (2014). Ching at al. (2012) investigated finite element simulation of mixed convection heat and mass transfer in a right triangular enclosure. His result depicts the increase of buoyancy ratio enhances the heat and mass transfer rate for all values of Richardson number and for each direction of the sliding wall motion.

Thermodynamic heat pump cycles or refrigeration cycles are the mathematical models for heat pump, air conditioning and refrigeration systems. A heat pump (heater) is a mechanical system that allows the transmission of heat from the source (at low temperature) to sink (at high temperature). Number of studies were investigated along with wall heater by the researcher such as MHD mixed convection in a lid-driven cavity with corner heater by Oztop et al (2011). Further investigation on MHD natural convection in an enclosure from two semi-circular heaters on the bottom wall was found by Oztop et al. (2012). Laminar mixed convection flow and heat transfer characteristics in a lid driven cavity with a circular cylinder studied by Khanafer and Aithel (2013). Also, Khanafer (2014) explored his research on comparison of flow and heat transfer characteristics in a lid-driven cavity between flexible and modified geometry of a heated bottom wall.

MHD studies are mostly focused on natural convection. For example, Temah (2008) investigated numerical simulation of double diffusive natural convection in rectangular enclosure in the presences of magnetic field and heat source. Rahman et al. (2011) studied MHD Mixed convection with joule heating effect in a lid-driven cavity with a heated semi-circular source using finite element technique.

Clearly, the combined heat and mass transfer on mixed convection lid-driven rectangular cavity containing heated obstacle have recognized great attention in the recent years. The main objective of this work is to present the effects of magnetic field and buoyancy ratio on laminar mixed convection problem foe a lid-driven cavity containing trapezoidal heated block and two semi-circular wall heaters on the heat and flow characteristics at two semi-circular wall heaters.

## 2. Physical Model and mathematical formulation:

The considered two-dimensional model is illustrated in Fig. 1 with boundary conditions. It is a rectangular enclosure with trapezoidal heated block in the middle of the enclosure and two semi-circular wall heaters at the bottom wall of the enclosure. As illustrates schematically, the top wall moves with constant velocity in the positive $x$-direction and the upper and bottom walls are kept adiabatic except two semi-circular wall heaters. Both vertical walls are subjected to low concentrated $C_L$ and low temperature $T_H$, the obstacle area and two semi-circular wall heaters are subjected to high concentrated $C_H$ and high temperature $T_H$. The heat transfer and fluid flow will be illustrated for commonly used fluid with $Pr = 7$. Hartmann number are considered $50 - 150$ to analyze the effect of magnetic field.

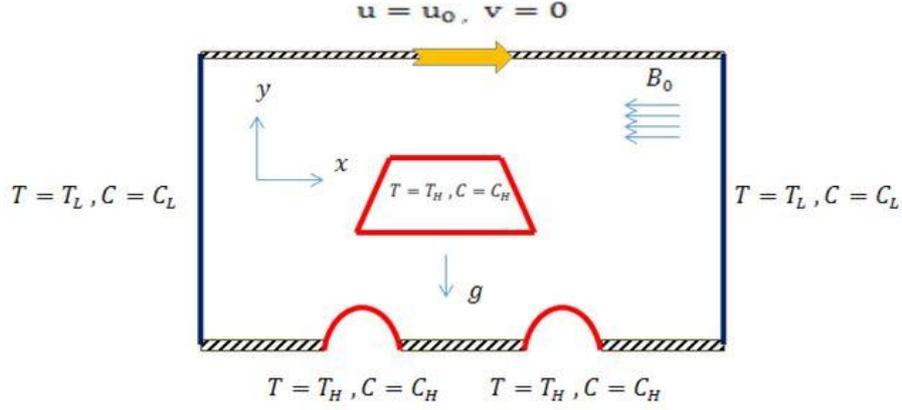

Fig. 1. Schematic diagram of the cavity with trapezoidal obstacle and two semi-circular wall heaters

Base on the model, two dimensional, laminar, steady equations are written by considering a uniform applied magnetic field. We assume that Boussinesq approximation is valid and the uniform magnetic field $B_0$ is also applied to the fluid in the direction parallel to at opposite of positive $x$ −direction. Thus, using the coordinate system shown in Fig.1, the governing equations can be written in dimensional form

$$\frac{\partial u}{\partial x} + \frac{\partial v}{\partial y} = 0 \tag{1}$$

$$\rho \left(u \frac{\partial u}{\partial x} + v \frac{\partial u}{\partial y}\right) = -\frac{\partial p}{\partial x} + \mu \left(\frac{\partial^2 u}{\partial x^2} + \frac{\partial^2 u}{\partial y^2}\right) \tag{2}$$

$$\rho \left(u \frac{\partial v}{\partial x} + v \frac{\partial v}{\partial y}\right) = -\frac{\partial p}{\partial y} + \mu \left(\frac{\partial^2 v}{\partial x^2} + \frac{\partial^2 v}{\partial y^2}\right) + \rho g \beta_T (T - T_L) + \rho g \beta_C (C - C_L) - \sigma B_0^2 v \tag{3}$$

$$\left(u \frac{\partial T}{\partial x} + v \frac{\partial T}{\partial y}\right) = \alpha \left(\frac{\partial^2 T}{\partial x^2} + \frac{\partial^2 T}{\partial y^2}\right) \tag{4}$$

$$\left(u \frac{\partial C}{\partial x} + v \frac{\partial C}{\partial y}\right) = D \left(\frac{\partial^2 C}{\partial x^2} + \frac{\partial^2 C}{\partial y^2}\right) \tag{5}$$

We introduce now the following dimensionless variables

$$X = \frac{x}{L}, Y = \frac{y}{L}, U = \frac{u}{u_0}, V = \frac{v}{u_0}, P = \frac{(p+\rho g y)}{\rho u_0^2}, \theta = \frac{T-T_L}{T_H-T_L}, C = \frac{C-C_L}{C_H-C_L}, Le = \frac{\alpha}{D}, Re = \frac{L u_0}{\nu},$$

$$\alpha = \frac{k}{\rho C_p}, \nu = \frac{\mu}{\rho} \tag{6}$$

Equation (1) to (5) become

$$\frac{\partial U}{\partial X} + \frac{\partial V}{\partial Y} = 0 \tag{7}$$

$$U \frac{\partial U}{\partial X} + V \frac{\partial U}{\partial Y} = -\frac{\partial P}{\partial X} + \frac{1}{Re} \left(\frac{\partial^2 U}{\partial X^2} + \frac{\partial^2 U}{\partial Y^2}\right) \tag{8}$$

$$U \frac{\partial V}{\partial X} + V \frac{\partial V}{\partial Y} = -\frac{\partial P}{\partial Y} + \frac{1}{Re} \left(\frac{\partial^2 V}{\partial X^2} + \frac{\partial^2 V}{\partial Y^2}\right) + Ri(\theta + Br\ C) - Ha^2 \frac{1}{Re} V \tag{9}$$

$$U \frac{\partial \theta}{\partial X} + V \frac{\partial \theta}{\partial Y} = \frac{1}{PrRe} \left(\frac{\partial^2 \theta}{\partial X^2} + \frac{\partial^2 \theta}{\partial Y^2}\right) \tag{10}$$

$$U \frac{\partial C}{\partial X} + V \frac{\partial C}{\partial Y} = \frac{1}{LePrRe} \left(\frac{\partial^2 C}{\partial X^2} + \frac{\partial^2 C}{\partial Y^2}\right) \tag{11}$$

The variables have their usual sense in fluid mechanics and heat transfer as listed in the nomenclature. Three parameters that preside over this problem are Richardson number (Ri), Hartmann number (Ha) and Buoyancy ratio (Br) which are defined as

$$Ri = \frac{g\beta_T(T_H-T_L)L}{u_0^2}, Ha^2 = \frac{\sigma B_0^2 L^2}{\mu} \text{ and } Br = \frac{\beta_C(C_H-C_L)}{\beta_T(T_H-T_L)}.$$

The corresponding boundary conditions for the above problem are given by:

at the left and right walls: $U = 0, V = 0, \theta = 0, C = 0 \quad \forall \quad 0 < Y < 1$ (12)

at the upper wall: $U = 1, V = 0, \frac{\partial \theta}{\partial y} = 0 \quad \forall \quad Y = 1$ (13)

on the middle obstacle and semi-circular wall heater:

$U = 0, V = 0, \theta = 1, C = 1 \quad \forall \quad 0 \leq X \leq 1, 0 \leq Y \leq 1$ (14)

At the bottom wall without semi-circular wall heater:

$U = 0, V = 0, \frac{\partial \theta}{\partial Y} = 0 \quad \forall \quad Y = 0$ (15)

The local Nusselt number is calculated by the following expression:

$Nu = -\frac{\partial \theta}{\partial n} = \sqrt{\left(\frac{\partial \theta}{\partial X}\right)^2 + \left(\frac{\partial \theta}{\partial Y}\right)^2}$. The average Nusselt number at the heated surface and semi-circular wall heater based on the dimensionless quantities may be expressed as $Nu_{av} = \frac{1}{L_S}\int_0^{L_S} Nu \, ds$ and $Nu_{av} = \frac{1}{L}\int_0^L Nu \, ds$ respectively where $L_S$ and $L$ denotes the arc length of the semi-circular wall heater and length of the heated wall respectively.

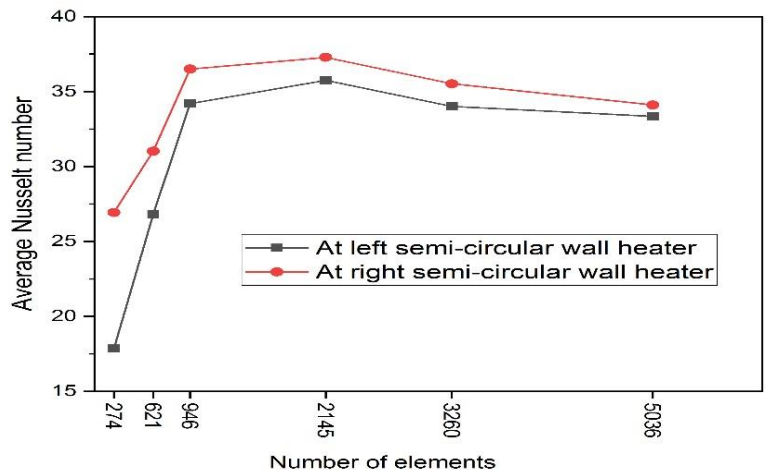

Fig. 2. Grid independency study: average Nusselt number at two semi-circular wall heaters for different grid elements while $Ri = 5, Br = 20$ and $Ha = 50$

## 3. Numerical solution:

### 3.1 Solution procedure

The computational procedure is similar to the works described by Ching et al (2012). The governing Eqs. $(8-11)$ along with the associated boundary conditions have been solved numerically by using the Galerkin weighted residual finite element method. The finite element method begins by the partition of the continuum area of interest into a number of simply shapes and sizes. Within each element, the dependent variables are approximated using interpolation functions. A non-uniform triangular mesh arrangement is implemented in the present investigation especially near the walls to capture the rapid changes in dependent variables. The velocity and thermal energy equations $(8) - (11)$ results in a set of non-linear coupled equations for which an iterative scheme is adopted. To ensure convergence of the numerical algorithm the following criteria is applied to all dependent variables over the solution domain $\sum |\psi_{ij}^n - \psi_{ij}^{n-1}| \leq 10^{-5}$, where $\psi$ represents a dependent variable $U, V, P, T$ and $C$; the indexes $i, j$ indicates a grid point and index $n$ is the current iteration at the grid level.

### 3.2. Grid refinement check:

In order to obtain grid independent solution, a grid refinement study is performed for a rectangular cavity with trapezoidal heated obstacle and two semi-circular wall heaters at $Ha = 50, Re = 100, Pr = 7, Ri = 5, Br = 20$ and $Le = 20$. In this investigation, six different non-uniform grids of triangular elements: 5036, 3260, 2145, 946, 621 and 274 are used. The values of average Nusselt number in presence of MHD at the two semi-circular wall heaters are used as a sensitivity measure of the accuracy of the solution and are selected at the monitoring variables for the grid independence study.

| Ri | Present | Y.C. Ching [21] |
|---|---|---|
| 0.01 | 30.258 | 32.386 |
| 0.1 | 27.687 | 28.653 |
| 1 | 12.323 | 12.231 |
| 10 | 11.029 | 11.5689 |

Table 1. Comparison of average Nusselt number between the present numerical solution and that of Ching et al (2012) at $Pr = 0.71, Re = 100$ and $Ha = 0$.

### 3.3 Code validation

The computational results are compared with the literature Ching et al (2012) for validate the present numerical code. The physical problem studied by Ching et al (2012) was a triangular cavity without MHD considered with fluid by finite element weighted residual method whose vertical wall was moved upward with a velocity and maintained at cooled condition. The inclined wall was hot, whereas the bottom wall was under the adiabatic conditions. Average Nusselt number is calculated for different values Ri whereas Pr and Re were kept atfixed conditions as shown in the Table 1. Validation of the code was also done by comparing streamline and isotherm in Fig. 3 by Ching et al (2012). As seen from this figure the obtained result shows very good agreement.

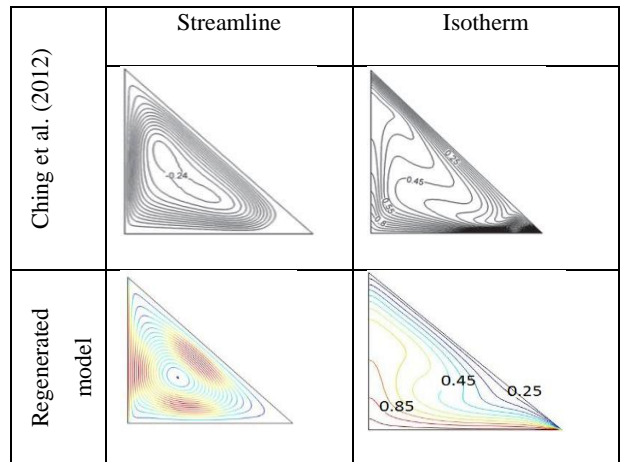

Fig. 3. Comparison of streamlines and isotherms with $Pr = 1, Br = 20, Le = 2$ and $Ri = 5$ with MHD

4. Results and discussion:

A numerical study has been performed in this work to investigate the effect of fluid flow and temperature in a lid-driven cavity in presence of magnetic field inside containing a trapezoidal heated block and two semi-circular wall heaters. Study is performed for different values of Hartmann number $Ha$, Richardson number $Ri$, Buoyancy ratio $Br$ and average Nusselt number $Nu_{av}$. Reynolds number and Lewis number are kept fixed as $Re = 100$ and $Le = 20$ respectively. For all cases Prandtl number is chosen $Pr = 7$ for water at 20

For the variation of $Ri$ and different values of $Ha$ with overall features, streamlines and isotherms are presented in Fig. 4 and Fig. 5 respectively to figure out the effects of the trapezoidal heated block and two semi-circular wall heaters on fluid flow and temperature distribution. At forced convection region in presence of MHD fluid flow is characterized by rotating vortex occupying the entire cavity and the density of fluid flow increased near the upper horizontal wall generated by the movement of upper lid. Also, a lonely vortex is pioneer at the left corner of the cavity near the left semi-circular wall heater. With the increase of $Ha$ from $50\ to\ 150$ the lonely vortex totally vanished and all the streamlines gathered on the top of the trapezoidal heated obstacle.

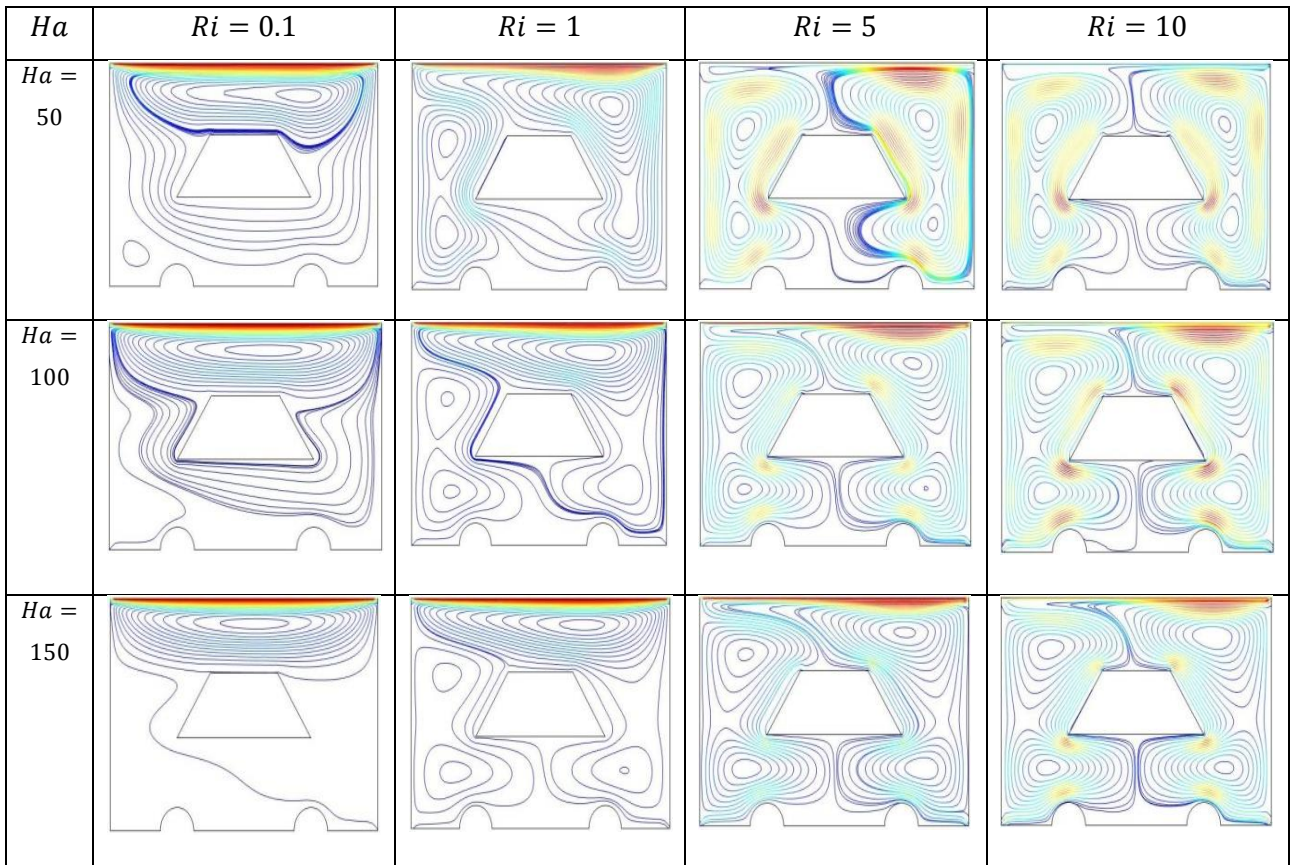

Fig. 4. Streamlines for different values of $Ri$ and $Ha$ for $Re = 100, Pr = 7, Le = 20, Br = 20$

At mixed convection region the core of the fluid flow divided into two parts positioned at the left and right side of the heated block. With the increase of $Ri$ and $Ha$ these divided parts become more prominent and symmetric seems like big mango seeds. Density of the streamlines near the moving lid decreases with the increases of $Ha$ and $Ri$. At mixed to natural convection region with the increase of $Ha$ the velocity field form two large bi-cellular vortex with four minor vortices. For highest value of $Ha$ and $Ri$ the appearance of streamlines become more symmetric due to the application of transverse magnetic fields which slow down the movement of the buoyancy-induced flow within the cavity.

The corresponding effect on the temperature fields shown in Fig. 4 illustrates that isotherms are almost parallel to both vertical walls for the highest value of $Ha$ at the forced convection region, indicating that most of the heat transfer process is carried out by conduction. However, some deviation in the conduction dominated region isotherms lines are initiated near the left top surface of the cavity. With the increase of $Ha$ in the mixed convection region, the patterns of the isotherms become linear to nonlinear zigzag shape. As $Ri$ and $Ha$ increases, the nonlinearity in the isotherms become higher and plume formation is profound at the left and right side of the heated block and other isotherm lines tend to parallel to the vertical walls. Moreover, the formation of the thermal boundary layers near the two vertical walls are to be initiated for the lower value of $Ha$. This is owing to the dominating influence of the convective current in the cavity. Also, for highest value of $Ha$ and $Ri$ the isotherm lines are in symmetric form.

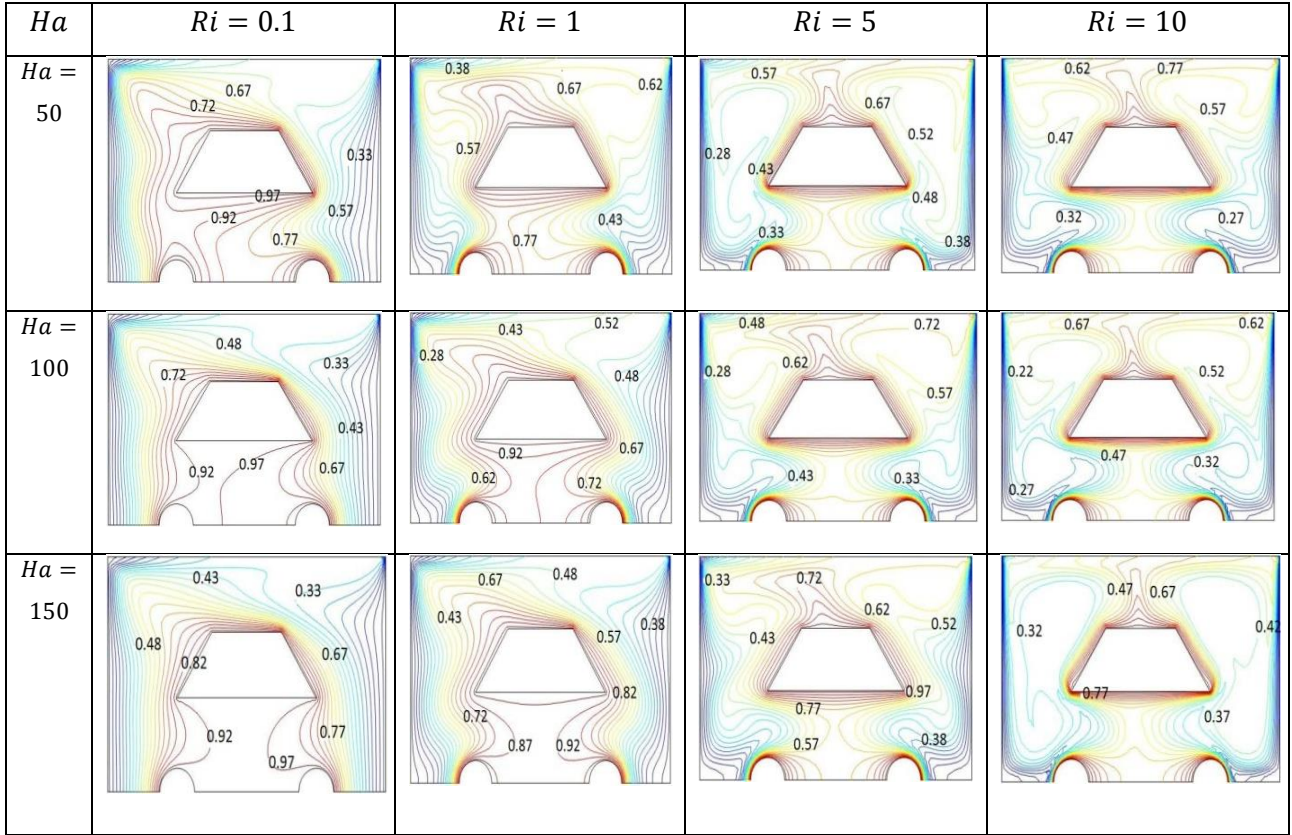

Fig. 5. Isotherms for different values of $Ri$ and $Ha$ for $Re = 100, Pr = 7, Le = 20, Br = 20$

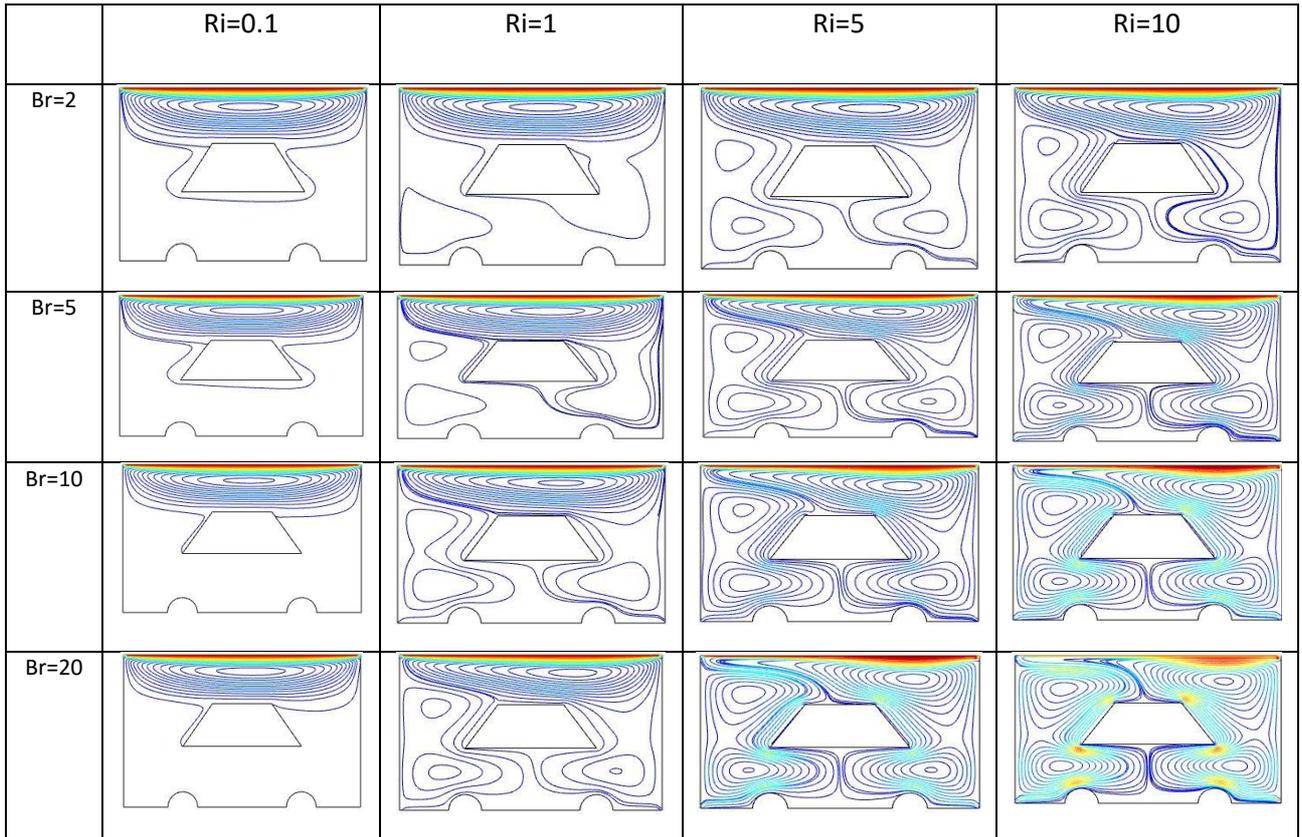

Fig. 6 Streamlines for different values of $Br$ and $Ri$ with $Ha = 150$ while $Re = 100, Pr = 7$

The effect of fluid flow for different values of $Br$ and $Ri$ with the trapezoidal obstacle and fixed $Ha = 150$, the overall features are illustrated in Fig. 5. At the forced convection region, all the streamlines are accumulated at the top of the heated block for each value of $Br$. At mixed convection region, with the increase of $Br$ a single vortex is formed. At highest value of $Ha$, the fluid flow experiences a Lorentz force due to the influence of the magnetic field. Further increase of $Ri$, streamline creates bi-cellular vortex at the left and right side of the block. The heated block divided the vortex into two parts and they appeared in symmetric condition. Also, the corresponding effect of temperature are shown in Fig. 6. Isotherm lines are almost parallel in forced and mixed convection region. With the increase of $Ri$ isotherm lines become nonlinear. At natural convection region the nonlinearity become higher and plume formation is profound at the left and right side of the heated block.

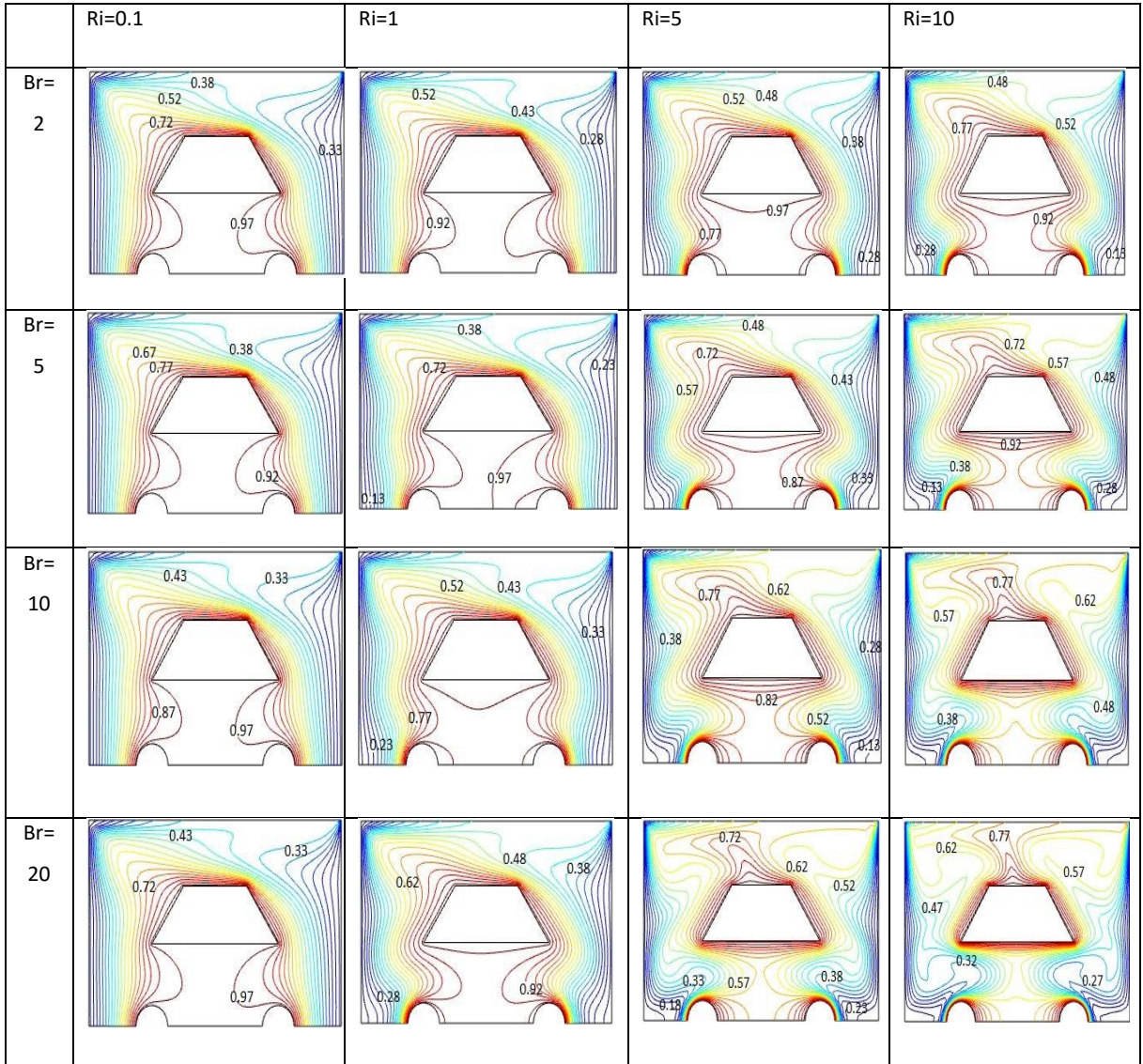

Fig. 7 Isotherms for different values of $Br$ and $Ri$ with $Ha = 150$ while $Re = 100, Pr = 7$

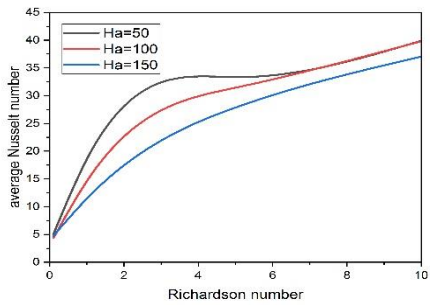                                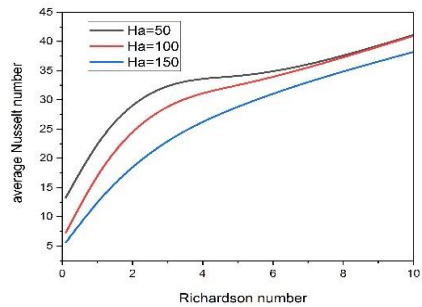

(a)                                                                                    (b)

Fig. 8. Average Nusselt number Vs Richardson number for different values of $Ha$ at (a) left semi-circular wall heater (b) right semi-circular wall heater

Figure 7 illustrates the average Nusselt number Vs Richardson number for different values of $Ha$ at left and right semi-circular wall heater. In presence of magnetic field it is quite normal to decrease heat transfer and with the increase of $Ha$ from 0 to 150 this procedure continues and it is observed that, for highest value of $Ha$ heat transfer rate is same also the curves are smother in both the cases. In the region $0.1 \leq Ri \leq 1$, the rate of heat transfer decreases slightly but it decreases abruptly in the region $1 \leq Ri \leq 10$. Also in forced convection region the average Nusselt number coincide for each values of $Ha$ at left semi-circular wall heater which is decreases slowly with the increase of $Ha$.

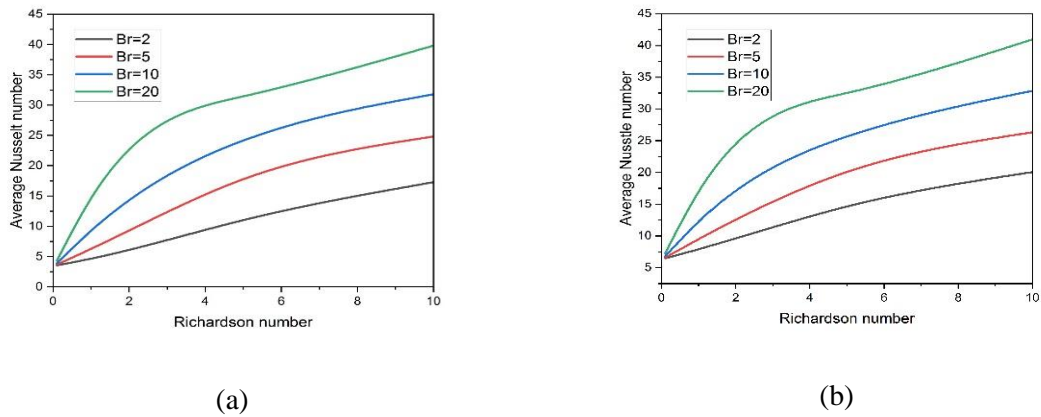

(a)　　　　　　　　　　　　　　　　(b)

Fig. 9. Average Nusselt number Vs Richardson number for different values of $Br$ at (a) left semi-circular wall heater (b) right semi-circular wall heater

Figure 8 illustrates the effect of $Br$ in presence of MHD with trapezoidal obstacle along the left and right semi-circular wall heater. This figure exhibits that, for each values of $Br$, heat transfer rate increases with the increase of $Ri$. Foe highest value of $Br$ heat transfer rate is high along the both left and right semi-circular wall heater. Average Nusselt number $Nu_{av}$ increases slightly in the region $0.1 \leq Ri \leq 1$ but increases abruptly in the region $1 \leq Ri \leq 10$. For lowest value of $Br$, the variation of $Nu_{av}$ is nominal in forced to mixed and mixed to natural convection region.

5. Conclusions:
The study has been analyzed with the numerical modeling of mixed convection of MHD flow in a rectangular lid-driven cavity containing trapezoidal heated obstacle in the middle and two semi-circular wall heaters at the bottom wall. The governing parameters that affect the flow and heat transfer characteristics are Hartmann number, Buoyancy ratio, Richardson number and two semi-circular wall heaters. In view of the obtained results, following findings can be drawn to make a summary as
- Magnetic field plays an important role on the flow pattern and temperature. They become weak for increasing the value of Hartmann number. In presence of magnetic field $Nu_{av}$ decreases drastically. Increasing value of Ha decreases $Nu_{av}$ gradually in the region $0.1 \leq Ri \leq 1$ but decreases $Nu_{av}$ abruptly in the region $1 \leq Ri \leq 10$.
- With the increase of $Br$ heat transfer rate increases. $Nu_{av}$ increases slightly in the region $0.1 \leq Ri \leq 1$ and it increases promptly in the region $1 \leq Ri \leq 10$ with    the increase of $Br$. Foe lower value of $Br$ heat transfer rate nominal in forced to mixed and mixed to natural convection region.
- Heat transfer rate is high along right semi-circular wall heater than the left semi-circular wall heater for increasing both parameter Hartmann number and Buoyancy ratio.
- There is no significant effect of Lewis number on fluid flow and temperature.